\documentclass[11pt]{article}

\usepackage{proof}
\usepackage{latexsym}

\newtheorem{theorem}{Theorem}[section]

\newtheorem{definition}[theorem]{Definition}
\newtheorem{proposition}[theorem]{Proposition}
\newtheorem{lemma}[theorem]{Lemma}

\title{
Iterating the recursively Mahlo operations
}
\author{
Toshiyasu Arai
\\
Graduate School of Science\\
Chiba University\\
1-33, Yayoi-cho, Inage-ku,
Chiba, 263-8522, JAPAN}
\date{}
\begin{document}
\maketitle

\begin{abstract}
In this paper we address a problem:
How far can we iterate lower recursively Mahlo operations in higher reflecting universes?
Or formally: 
How much can lower recursively Mahlo operations be iterated in set theories for 
higher reflecting universes?

It turns out that in $\Pi_{N}$-reflecting universes 
the lowest recursively Mahlo operation can be iterated along 
 towers of $\Sigma_{1}$-exponential orderings of height $N-3$,
and that all we can do is such iterations.
Namely
 the set theory for $\Pi_{N}$-reflecting universes 
is proof-theoretically reducible to 
 iterations of the operation along 
 such a tower.
\end{abstract}

For set-theoretic formulas $\varphi$, 
\[
P\models\varphi:\Leftrightarrow (P,\in)\models\varphi
.\]

In what follows, let L denote a transitive set, which is a universe in discourse.
$P, Q,\ldots$ denotes transitive sets in $\mbox{L}\cup\{\mbox{L}\}$ such that $\omega\in P$.

Let ${\cal X}$ be a first-order class of transitive sets.
This means that there exists a first-order sentence $\varphi$ such that $P\in {\cal X} \Leftrightarrow P \models \varphi$.
Then a set theory T is said to prove $\mbox{L}\in{\cal X}$ iff $\mbox{T}\vdash\varphi$.

A $\Pi_{i}$-recursively Mahlo operation for $2\leq i<\omega$, is then defined through a universal $\Pi_{i}$-formula $\Pi_{i}(a)$: 
\begin{eqnarray*}
P\in M_{i}({\cal X}) & :\Leftrightarrow & \forall b\in P[P\models\Pi_{i}(b) \to \exists Q\in {\cal X}\cap P(b\in Q\models\Pi_{i}(b))]
\\
&&\mbox{(read:} P \mbox{ is } \Pi_{i}\mbox{-reflecting on } {\cal X}\mbox{.)}
\end{eqnarray*}

Its iteration is defined by transfinite recursion on ordinals $\beta$:
\[M_{i}^{\beta}:=
\bigcap\{M_{i}(M_{i}^{\nu}):
 \nu<\beta\}.\]
Observe that $M_{i}({\cal X})$ is $\Pi_{i+1}$, i.e., there exists a $\Pi_{i+1}$-sentence $m_{i}({\cal X})$
such that $P\in M_{i}({\cal X})$ iff $P\models m_{i}({\cal X})$ for any transitive (and admissible) set $P$.

A transitive set $P$ is said to be $\Pi_{i}${\it -reflecting\/} 
if $P\in M_{i}=M_{i}^{1}$.

Let us denote
\[
{\cal X}\prec_{i}{\cal Y}  :\Leftrightarrow  {\cal Y}\subseteq M_{i}({\cal X})\mbox{, i.e., } \forall P\in{\cal Y}(P\in M_{i}({\cal X}))
.\]

$P\in M_{i+1}$ is much stronger than $P\in M_{i}$: Assume $P\in M_{i+1}$ and 
$P\models\Pi_{i}(b)$ for $b\in P$.
Then $P\in M_{i}$ and $P\models m_{i}\land\Pi_{i}(b)$ for the $\Pi_{i+1}$-sentence $m_{i}$
such that $P\in M_{i}$ iff $P\models m_{i}$.
Hence there exists a $Q\in P$ such that $Q\models m_{i}\land\Pi_{i}(b)$, i.e.,
$Q\in M_{i} \,\&\, Q\models\Pi_{i}(b)$.
This means $P\in M_{i}^{2}=M_{i}(M_{i})$, i.e., $M_{i}\prec_{i}M_{i+1}$.
Moreover
$P\in M_{i}^{\triangle}$, i.e., $P\in\bigcap\{M_{i}^{\beta}:\beta\in ord(P)\}$, 
$M_{i}^{\triangle}\prec_{i}M_{i+1}$, and so on.

In particular a set theory KP$\Pi_{i+1}$ for universes in $M_{i+1}$ proves the consistency of a set theory 
for universes in $M_{i}^{\triangle}$.

In this paper we address a problem:
How far can we iterate lower recursively Mahlo operations in higher reflecting universes?
Or formally: 
How much can lower recursively Mahlo operations be iterated in set theories for 
higher reflecting universes?
Specifically:
What kind of iterations of the lowest operations $M_{2}$ do we need to obtain equiconsistent
theories for set theories for higher reflecting universes?

\section{Iterations of the operation $M_{i}$ in $\Pi_{i+1}$-reflectings}

In this section we see that 
iterations of the operation $M_{i}$ along $\Sigma_{1}$-relations on $\omega$
are too short to resolve $\Pi_{i+1}$-reflecting universes 
provided that the $\Sigma_{1}$-relations are provably wellfounded in KP$\Pi_{i+1}$.

\begin{definition}\label{df:LE0}
\begin{enumerate}
\item
\label{df:LE1}
{\rm KP}$\ell$ {\rm denotes a set theory for limits of admissibles.}
 {\rm KP}$\Pi_{N}$ {\rm denotes a set theory for universes in} $M_{N}$.
 
 \item\label{df:LE2}
{\rm For a definable relation} $\prec$
{\rm and set-theoretic universe} $P$ {\rm (admissibility suffices)}
{\rm let}
\[
P\in M_{i}(a;\prec) :\Leftrightarrow P\in\bigcap\{ M_{i}(M_{i}(b;\prec)): b \prec^{P} a\},
\]
{\rm where} $b\prec^{P}a:\Leftrightarrow P\models b\prec a$.

{\rm Note that} $M_{i}(a;\prec)$ {\rm is a} $\Pi_{i+1}${\rm -class for (set-theoretic)} $\Sigma_{i+1}$ $\prec$.

\item
{\rm We say that a theory T is} proof-theoretically reducible {\rm to another theory S
if T is a} $\Pi^{1}_{1}$ {\rm (on} $\omega${\rm )-conservative extension of S, and the fact is provable
in a weak arithmetic, e.g., the elementary recursive arithmetic EA.}

\item
{\rm For a relation} $\prec$ {\rm on} $\omega${\rm ,} $TI(a, \prec)$ {\rm denotes the transfinite induction schema up to} $a\in\omega${\rm :}
\[
\{ \forall x\in\omega[\forall y\prec x \varphi(y) \to \varphi(x)] \to \forall x\prec a \varphi(x):
 \varphi \mbox{ {\rm is a set-theoretic formula}}\}
\]
{\rm and} $TI(a, \prec,\Pi_{n})$ {\rm its restriction to} $\Pi_{n}${\rm -formulas} $\varphi$.

{\rm Using a universal} $\Pi_{n}${\rm -formula,} $TI(a, \prec,\Pi_{n})$ {\rm is equivalent to a single}
$\Pi_{n+2}${\rm -formula.}

\item
{\rm A relation} $\prec$ {\rm on} $\omega$ {\rm is said to be} almost wellfounded {\rm in KP}$\ell$
{\rm if KP}$\ell$ {\rm proves the transfinite induction schema}
$TI(a,\prec)$
{\rm up to} each $a\in\omega$.
 \end{enumerate}
\end{definition}

It is easy to see the following lemma using the fact that $M_{i}(a;\prec)$ is $\Pi_{i+1}$.

\begin{lemma}\label{lem:resolution-1}
Let $\prec$ be a $\Sigma_{1}$ relation on $\omega$.
Then {\rm KP}$\Pi_{i+1}\, (i\geq 2)$ proves
\[
\forall a\in\omega[TI(a, \prec,\Pi_{i+1}) \to \mbox{{\rm L}}\in M_{i}(a;\prec)]
.\]
A fortiori
{\rm KP}$\Pi_{i+1}$ proves
$\forall a\in\omega[TI(a, \prec,\Pi_{i+1}) \to \mbox{{\rm L}}\in M_{2}(M_{i}(a;\prec))]$.

In other words, {\rm KP}$\ell$ proves
$P\in M_{i+1} \to \forall a\in\omega[TI(a, \prec^{P},\Pi_{i+1}^{P}) \to P\in M_{i}(a;\prec)]$.
\end{lemma}

Therefore $\forall a\in\omega[\mbox{{\rm L}}\in M_{i}(a;\prec)]$ is too weak to reduce
KP$\Pi_{i+1}$ proof-theoretically for any $\Sigma_{1}$ relation $ \prec$ on $\omega$, for example
$\mbox{KP}\Pi_{i+1}\vdash\mbox{CON}(\forall a\in\omega[\mbox{{\rm L}}\in M_{i}(a;\prec)])$
if $\forall a\in\omega[TI(a,\prec)]$ is provable in KP$\Pi_{i+1}$.

Nonetheless $\Pi_{i+1}$-reflecting universes can be approximated by iterations of the operation $M_{i}$
along well founded $\Sigma_{1}$ relations on $\omega$.

\begin{theorem}\label{th:resolution0}
For each $i\, (2\leq i<\omega)$ there exists a $\Sigma_{1}$ almost wellfounded relation $\lhd_{i}$ in {\rm KP}$\ell$
such that
{\rm KP}$\Pi_{i+1}$ is proof-theoretically reducible to
the theory 
\[\mbox{{\rm KP}}\ell+
\{
\mbox{{\rm L}}\in M_{i}(a;\lhd_{i}): 
a\in\omega\}
.\]
\end{theorem}

Theorem \ref{th:resolution0} follows from Lemma \ref{lem:i+1toi} and Theorem  \ref{th:resolution2} below.

The case $i=2$ means that the set theory KP$\Pi_{3}$ for $\Pi_{3}$-reflecting universes 
can be resolved by iterations of the recursively
Mahlo operations $M_{2}$.
\\ \smallskip \\
{\bf Remark}.
Although KP$\ell$ is weaker than KP$\Pi_{i+1}$, KP$\Pi_{i+1}$ does not prove the soundness of KP$\ell$:
Let Fund denote the axiom schema for Foundation.
Then for a $\varphi\in\Sigma_{i+2}$ and a standard provability predicate 
$\mbox{Pr}_{\mbox{{\footnotesize Fund}}}$ of Fund
\[
\mbox{KP}\Pi_{i+1}\not\vdash \forall n\in\omega[
\mbox{Pr}_{\mbox{{\footnotesize Fund}}}(\lceil\varphi(\dot{n})\rceil)
\to \varphi(n)]
\]
since $\mbox{KP}\Pi_{i+1}\setminus \mbox{Fund}\subseteq\Pi_{i+2}\, (i\geq 2)$.

Hence even if $\mbox{KP}\Pi_{i+1}\vdash 
\forall a\in\omega[\mbox{Pr}_{\mbox{{\footnotesize KP}}\ell}(\lceil TI(\dot{a}, \lhd_{i},\Pi_{i+1})\rceil)]$,
this does not imply $\mbox{KP}\Pi_{i+1}\vdash \forall a\in\omega\, TI(a, \lhd_{i},\Pi_{i+1})$.

\section{$\Pi_{3}$-reflecting on $\Pi_{3}$-reflectings}

Our goal is to approximate $\Pi_{i+1}$-reflecting universes by iterations of the lowest recursively
Mahlo operations $M_{2}$.
Let us consider first the simplest case: 
$\Pi_{3}$-reflecting universes on $\Pi_{3}$-reflectings, $M_{3}^{2}=M_{3}(M_{3})$.
Universes in $M_{3}^{2}$ are seen to be resolved in terms of iterations of the operation $M_{2}$ along 
a {\it lexicographic ordering\/} on pairs.

\begin{definition}\label{df:LEx}
\begin{enumerate}

\item\label{df:LE3}
{\rm For a} $\Sigma_{1}$ {\rm relation} $\prec$ {\rm on} $\omega${\rm ,}
$W(\prec)$ {\rm  denotes the} wellfounded part {\rm of} $\prec$:
\[
a\in W(\prec) :\Leftrightarrow \forall f\in{}^{\omega}\omega\exists n\in\omega[f(0)=a \to f(n+1)\not\prec f(n)]
.\]
$W(\prec)$ {\rm is} $\Pi_{1}$.

{\rm Note that} $W(\prec^{Q})$ {\rm is a} set {\rm in limits of admissibles} $P$ 
{\rm for any transitive set} $Q\in P$.

\item\label{df:LE4}
{\rm For two transitive relations} $<_{1}, <_{0}$ {\rm on} $\omega$,
 $<_{L}:\equiv L(<_{1},<_{0})$ {\rm denotes the lexicographic ordering:}
\[
\langle n_{1},n_{0}\rangle<_{L}\langle m_{1},m_{0}\rangle :\Leftrightarrow n_{1}<_{1}m_{1} \mbox{ {\rm or} }
( n_{1}=m_{1} \,\&\, n_{0}<_{0}m_{0})
.\]
$L(<_{1},<_{0})$ {\rm is} $\Sigma_{1}$ {\rm if} $<_{1}$ {\rm and} $<_{0}$ {\rm are} $\Sigma_{1}$.

 $<_{LW}$ {\rm denotes the restriction of} $<_{L}$ 
 {\rm to the wellfounded part in the second component:}
\[
\langle n_{1},n_{0}\rangle<_{LW}\langle m_{1},m_{0}\rangle :\Leftrightarrow \langle n_{1},n_{0}\rangle<_{L}\langle m_{1},m_{0}\rangle 
\,\&\, n_{0},m_{0}\in W(<_{0})
.\]
$<_{LW}$ {\rm is} $\Delta_{2}$ {\rm if} $<_{1}$ {\rm and} $<_{0}$ {\rm are} $\Sigma_{1}$.
\end{enumerate}
\end{definition}

\begin{proposition}\label{prp:axbeta}
Let $P$ be a limit of admissibles and $<$ be a $\Sigma_{1}$ relation on $\omega$.
Suppose $P\models a\in W(<)$.
Then $a\in W^{P}(<^{Q})=W(<^{Q})$ and $Q\models TI(a, <)$ for any $Q\in P$,
where
\[
a\in W^{P}(<^{Q}) :\Leftrightarrow \forall f\in{}^{\omega}\omega\cap P\exists n\in\omega[f(0)=a \to f(n+1)\not<^{Q} f(n)]
.\]
\end{proposition}
{\bf Proof}.\hspace{2mm}
Since $<$ is $\Sigma_{1}$ and $Q\subseteq P$, we have $<^{Q}\subseteq <^{P}$.
Hence $a\in W^{P}(<^{P})\subseteq W^{R}(<^{Q})$ for any $R\subseteq P$.
Therefore $a\in W^{P}(<^{Q})=W^{Q^{+}}(<^{Q})=W(<^{Q})$ for the set $<^{Q}$ in $P$, and 
the next admissible $Q^{+}\in P$ above $Q$.
This yields the transfinite induction
schema $TI(a,<^{Q})$ up to $a$.
\hspace*{\fill} $\Box$

KP$\Pi_{3}(\Pi_{3})$ denotes a set theory for universes in
$M_{3}(M_{3})$.

\begin{lemma}\label{lem:m3m3}
Let $<_{1}, <_{0}$ be two $\Sigma_{1}$ transitive relations on $\omega$, and $<_{LW}$ the restriction of the lexicographic ordering defined from these to 
the wellfounded part in the second component.

Then {\rm KP}$\Pi_{3}(\Pi_{3})$ proves
\[
\forall a, \alpha\in\omega[TI(a, <_{1},\Pi_{3}) \to \mbox{{\rm L}}\in M_{2}(\langle a,\alpha\rangle;<_{LW})]
.\]
\end{lemma}
{\bf Proof}.\hspace{2mm}
Let $\mbox{L}\in M_{3}(M_{3})$.
By transfinite induction on $a$ along $<_{1}$ we show
\[
\forall \alpha\in\omega[\mbox{{\rm L}}\in M_{2}(\langle a,\alpha\rangle;<_{LW})]
\]
where
\[
P\in M_{2}(\langle a,\alpha\rangle;<_{LW}) \Leftrightarrow 
P\in\bigcap\{ M_{2}(M_{2}(\langle b,\beta\rangle; <_{LW})): 
 \langle b,\beta\rangle<_{LW}^{P}\langle a,\alpha\rangle\}
 \]
 and
 \[
  \langle b,\beta\rangle<_{LW}^{P}\langle a,\alpha\rangle \Leftrightarrow  
  \langle b,\beta\rangle<_{L}^{P}\langle a,\alpha\rangle \,\&\, P\models \alpha,\beta\in W(<_{0})
  .\]
  
Suppose that $\forall b<_{1}a\forall \beta\in \omega[\mbox{{\rm L}}\in M_{2}(\langle b,\beta\rangle;<_{LW})]$, 
and $\langle b,\beta\rangle<_{LW}\langle a,\alpha\rangle$.
We show
$\mbox{{\rm L}}\in M_{2}(M_{2}(\langle b,\beta\rangle;<_{LW}))$.

IH yields the case $b<_{1}a$. Assume $b=a$ and $\beta<_{0}\alpha\in W(<_{0})$.
Suppose a $\varphi\in\Pi_{2}$ holds in $\mbox{L}\in M_{3}(M_{3})$.
Pick a $Q\in \mbox{L}\cap M_{3}$ so that $Q\models\varphi$ and
$Q\in \bigcap\{ M_{2}(M_{2}(\langle b,\gamma\rangle; <_{LW})): 
Q\models b<_{1}a \land \gamma\in W(<_{0})\}$ by IH.

We claim that $Q\in M_{2}(\langle a,\beta\rangle;<_{LW})$.
By Proposition \ref{prp:axbeta} we have $Q\models TI(\beta,<_{0})$.
Hence we have
$Q\in M_{2}(\langle a,\beta\rangle;<_{LW})$
by transfinite induction on $\beta$.
\hspace*{\fill} $\Box$

\begin{theorem}\label{th:resolution1}

There exist $\Sigma_{1}$ transitive relations $<_{1}, <_{0}$ on $\omega$ such that
$<_{1}$ is almost wellfounded in {\rm KP}$\ell$, and
{\rm KP}$\Pi_{3}(\Pi_{3})$ is proof-theoretically reducible to
the theory 
\[
\mbox{{\rm KP}}\ell+
\{
\mbox{{\rm L}}\in \bigcap\{M_{2}(M_{2}(\langle a,\alpha\rangle; <_{LW})): 
\alpha\in W(<_{0})\}
: a\in\omega\}
\]
for  the restriction $<_{LW}$ of the lexicographic ordering $<_{L}=L(<_{1},<_{0})$ 
defined from these to 
the wellfounded part in the second components.
\end{theorem}

For a proof of Theorem \ref{th:resolution1}, see \cite{WienpiN}.

\section{$\Pi_{N}$-reflection}

As you expected, an {\it exponential\/} structure involves in resolving $\Pi_{N}$-reflecting universes L.

\begin{definition}\label{df:LEe}
{\rm Let} $<_{1}, <_{0}$ {\rm be two transitive relations on} $\omega$.
\begin{enumerate}
\item
{\rm The relation} $<_{E}=E(<_{1},<_{0})$ 
 {\rm  is on sequences} $\langle (n^{1}_{i}, n^{0}_{i}): i<\ell\rangle$ 
{\rm of pairs with} $<_{1}${\rm -decreasing first components} 
 {\rm (}$n^{1}_{i+1}<_{1}n^{1}_{i}${\rm ), and is defined by}
 \begin{eqnarray*}
 && \langle (n^{1}_{i}, n^{0}_{i}): i<\ell_{0}\rangle <_{E}
 \langle (m^{1}_{i}, m^{0}_{i}): i<\ell_{1}\rangle
 \mbox{ {\rm iff}}
 \\
 \mbox{{\rm either}} &&
 \\
 &&
 \exists k\forall i<k\forall j<2[n^{j}_{i}=m^{j}_{i} \,\&\, 
 (n^{1}_{k}, n^{0}_{k})<_{L}(m^{1}_{k},m^{0}_{k})]
 \\
 \mbox{ {\rm or}} && 
 \\
 &&
 \ell_{0}<\ell_{1} \,\&\, \forall i<\ell_{0}\forall j<2[n^{j}_{i}=m^{j}_{i}]
 \end{eqnarray*}
{\rm where} $<_{L}=L(<_{1},<_{0})$ {\rm in Definition \ref{df:LEx}.\ref{df:LE4}.}

{\rm Write}
$\sum_{i<\ell}\pi^{n^{1}_{i}}n^{0}_{i}$
{\rm for}
$\langle (n^{1}_{i}, n^{0}_{i}): i<\ell\rangle$.

\item
{\rm Let} $dom(<_{E})$ {\rm denote the domain of the relation} $<_{E}${\rm :}
\[
dom(<_{E}):=\{\sum_{i<\ell}\pi^{n^{1}_{i}}n^{0}_{i} : \forall i<\ell\dot{-}1(n^{1}_{i+1}<_{1}n^{1}_{i}) \,\&\,
n^{1}_{i}, n^{0}_{i},\ell\in\omega\}
.\]
\item
 $<_{EW}$ {\rm denotes the restriction of} $<_{E}$ 
 {\rm to the wellfounded part in the second components:}
\begin{eqnarray*}
&&
\alpha=\sum_{i<\ell_{0}}\pi^{n^{1}_{i}}n^{0}_{i} <_{EW} \sum_{i<\ell_{1}}\pi^{m^{1}_{i}}m^{0}_{i} =\beta
 \mbox{ {\rm iff}}
 \\
 &&
\alpha<_{E} \beta
\,\&\, \{n^{0}_{i}:i<\ell_{0}\}\cup\{m^{0}_{i}:i<\ell_{1}\}\subseteq W(<_{0})
.\end{eqnarray*}
\end{enumerate}
\end{definition}

\begin{lemma}\label{lem:i+1toi}
Let $<_{1}, <_{0}$ be two transitive relations on $\omega$, $<_{1}$ is $\Delta_{2}$, 
$<_{0}$ is $\Sigma_{1}$, 
and $<_{EW}$ the restriction of the exponential ordering defined from these to 
the wellfounded part in the second components.
Then {\rm KP}$\ell$ proves for each $i\geq 2$
\[
\forall P\in\mbox{{\rm L}}\cup\{\mbox{{\rm L}}\}\forall a\in\omega\forall \alpha<^{P}a[
P\in M_{i+1}(M_{i+1}(a;<_{1}))  \to P\in M_{i}(\alpha;<_{EW})]
\]
where for $\alpha=\sum_{i<\ell}\pi^{n^{1}_{i}}n^{0}_{i}\in dom(<^{P}_{E})$,
$\alpha<^{P}a :\Leftrightarrow n^{1}_{0}<^{P}_{1}a$.
\end{lemma}
{\bf Proof}.\hspace{2mm}
We show for any $a\in\omega$ and any $\beta\in dom(<_{EW}^{P}\uparrow a)$
\[
P\in M_{i+1}(M_{i+1}(a;<_{1})) \,\&\, P\in M_{i}(\beta;<_{EW})  
\to \forall\alpha<^{P}a\{P\in M_{i}(\beta+\alpha;<_{EW})\}
\]
by main induction on $P\in\mbox{{\rm L}}\cup\{\mbox{{\rm L}}\}$ with respect to the relation $\in$,
where for $\beta=\sum_{i<\ell_{1}}\pi^{m^{1}_{i}}m^{0}_{i}$ and 
$\alpha=\sum_{i<\ell_{0}}\pi^{n^{1}_{i}}n^{0}_{i}$,
\[
\beta\in dom(<_{EW}^{P}\uparrow a) :\Leftrightarrow 
\beta\in dom(<_{EW}^{P}) \,\&\, (\ell_{1}>0 \to a\leq_{1}^{P}m^{1}_{\ell_{1}-1})
\]
and
$\beta+\alpha=\sum_{i<\ell_{1}}\pi^{m^{1}_{i}}m^{0}_{i}+\sum_{i<\ell_{0}}\pi^{n^{1}_{i}}n^{0}_{i}$.

Suppose $\beta\in dom(<_{EW}^{P}\uparrow a)$, 
$P\in M_{i+1}(M_{i+1}(a;<_{1}))$ and  
$P\in M_{i}(\beta;<_{EW})$.
Pick an $\alpha=\pi^{b}x+\alpha_{0}\in dom(<_{EW}^{P})$ so that $\alpha_{0}<^{P}b<_{1}^{P}a$ and 
$x\in W^{P}(<_{0}^{P})$.
We show $P\in M_{i}(\beta+\alpha;<_{EW})$.
It suffices to show $P\in M_{i}(M_{i}(\beta+\gamma;<_{EW}))$ for any $\gamma<_{EW}^{P}\alpha$
by $P\in M_{i}(\beta;<_{EW})$.

If $\gamma$ is the empty sequence, then 
$P\in M_{i}(M_{i}(\beta;<_{EW}))$ follows from $P\in M_{i}(\beta;<_{EW})$, which is $\Pi_{i+1}$,
and $P\in M_{i+1}(M_{i+1}(a;<_{1}))\subseteq M_{i+1}$.

Let $\gamma=\pi^{c}y+\gamma_{0}$ with $\gamma_{0}<^{P}c\leq_{1}^{P}b$, 
and $P\models \theta$ for a $\theta\in\Pi_{i}$.
It suffices to find a $Q\in P$ so that $Q\in M_{i}(\beta+\gamma;<_{EW})$ and $Q\models\theta$.

First consider the case when $c<_{1}^{P}b$.
By $P\in M_{i+1}(M_{i+1}(a;<_{1}))$, pick a $Q\in P$ so that $Q\in M_{i+1}(a;<_{1})$, $Q\models\theta$,
$\beta\in dom(<_{EW}^{Q}\uparrow a)$, $Q\in M_{i}(\beta;<_{EW})$ and 
$dom(<_{EW}^{Q})\ni\gamma<^{Q}b<_{1}^{Q}a$.

Then $Q\in M_{i+1}(a;<_{1})\subseteq M_{i+1}(M_{i+1}(b;<_{1}))$, and hence MIH yields
$Q\in M_{i}(\beta+\gamma;<_{EW})$.

Thus we have shown $P\in \bigcap\{M_{i}(\beta+\delta;<_{EW}): \delta<^{P}b\}$, which is $\Pi_{i+1}$,
and hence
\begin{equation}\label{eq:MIH}
P\in M_{i}(M_{i+1}(a;<_{1})\cap \bigcap\{ M_{i}(\beta+\delta;<_{EW}): \delta<b\})
\end{equation}

Second consider the case when $c=b$.

We can find a $Q\in P$ so that $Q\in M_{i+1}(a;<_{1})$, $Q\models\theta$,
$\beta\in dom(<_{EW}^{Q}\uparrow a)$, $Q\in \bigcap\{M_{i}(\beta+\delta;<_{EW}): \delta<^{Q}b\}$ 
by (\ref{eq:MIH}) 
and 
$dom(<_{EW}^{Q})\ni\gamma \,\&\, b<_{1}^{Q}a$. 
We have $x\in W^{P}(<_{0}^{P})\subseteq W(<_{0}^{Q})$ by Proposition \ref{prp:axbeta}.

Therefore it suffices to show
\begin{eqnarray*}
&& \forall x\in W(<_{0}^{Q})\forall b\in\omega\forall \beta\in dom(<_{EW}^{Q}\uparrow b)
[Q\in P \,\&\,
Q\in M_{i+1}(M_{i+1}(b;<_{1})) \,\&\, \\
&& Q\in \bigcap\{M_{i}(\beta+\delta;<_{EW}): \delta<^{Q}b\}
\\
&&
\Longrightarrow  \forall\gamma_{0}<^{Q}b\{Q\in M_{i}(\beta+\pi^{b}x+\gamma_{0};<_{EW})\}]
\end{eqnarray*}
by subsidiary induction on $x\in W(<_{0}^{Q})$.

First assume $\beta+\pi^{b}y+\delta_{0}<_{EW}^{Q}\beta+\pi^{b}x+\gamma_{0}$ with $y<_{0}^{Q}x$.
SIH yields $Q\in M_{i}(\beta+\pi^{b}y+\delta_{0};<_{EW})$, and this implies
$Q\in M_{i}(M_{i}(\beta+\pi^{b}y+\delta_{0};<_{EW}))$ by $Q\in M_{i+1}$.

Therefore we have shown $Q\in M_{i}(\beta+\pi^{b}x;<_{EW})$ with $\gamma_{0}=0$.
Now let $\gamma_{0}=\pi^{c}y+\gamma_{1}$ with $c<_{1}^{Q}b$.
We have $\beta+\pi^{b}x\in dom(<_{EW}^{Q}\uparrow c)$, 
$Q\in M_{i+1}(M_{i+1}(b;<_{1})) \,\&\, Q\in M_{i}(\beta+\pi^{b}x;<_{EW})$ and $Q\in P$.
Hence MIH yields $Q\in M_{i}(\beta+\pi^{b}x+\gamma_{0};<_{EW})$ for $\gamma_{0}<^{Q}b$.
\hspace*{\fill} $\Box$

\begin{definition}\label{df:tower}
{\rm Let} $<_{i}\, (2\leq i\leq N-1)$ {\rm be} $\Sigma_{1}$ {\rm relations on} $\omega$.
{\rm Define a} tower {\rm relation} $<_{T}$ {\rm from these as follows.}

{\rm Define inductively relations} $<_{E_{i}}\, (2\leq i\leq N-1)$.
 \begin{enumerate}
 \item
 $<_{E_{N-1}}:\equiv <_{N-1}$.
 \item
 $<_{E_{i}}:\equiv E(<_{E_{i+1}},<_{i})$ {\rm for} $2\leq i\leq N-2${\rm , cf. Definition \ref{df:LEe}.}
 \end{enumerate}
{\rm Then let}
\[
<_{T}:\equiv <_{E_{2}}
.\]

$<_{TW}$ {\rm denotes the restriction of} $<_{T}$ {\rm to the wellfounded parts in the second components hereditarily.
Namely}
$<_{TW}=<_{E_{2}W}$ {\rm and} 
\begin{eqnarray*}
&& \sum_{n<\ell}\pi^{\alpha_{n}}x_{n}\in dom(<_{E_{i}W}) :\Leftrightarrow
\\
&&
\forall n<\ell\dot{-}1(\alpha_{n+1}<_{E_{i+1}W}\alpha_{n}) \,\&\, \forall n<\ell(x_{n}\in W(<_{i}))
\end{eqnarray*}
{\rm with} $<_{E_{N-1}W}=<_{N-1}$.

{\rm For} $a\in\omega$ {\rm and} $\alpha=\sum_{n<\ell}\pi^{\alpha_{n}}x_{n}\in dom(<_{T})${\rm , define inductively}
\[
\alpha<a
:\Leftrightarrow
\forall n<\ell(\alpha_{n}<a)
\]
{\rm with}
$\alpha_{n}<a :\Leftrightarrow \alpha_{n}<_{N-1}a$ {\rm for}
$\alpha_{n}\in\omega$.
\end{definition}

Lemmas \ref{lem:i+1toi} and \ref{lem:resolution-1} yield the following for the set theory KP$\Pi_{N}$ 
for universes in
$M_{N}$.

\begin{theorem}\label{th:tower}
Let $<_{i}\, (2\leq i\leq N-1<\omega)$ be $\Sigma_{1}$ transitive relations  on $\omega$.
Let $<_{TW}$ denote  the restriction of the tower $<_{T}$ of the exponential orderings $<_{E_{i}}$
defined from these to 
the wellfounded parts in the second components hereditarily.

Then
{\rm KP}$\Pi_{N}$ proves that
\[
\forall a\in\omega\forall \alpha<a [TI(a,<_{N-1},\Pi_{N}) \to \mbox{{\rm L}}\in M_{2}(\alpha;<_{TW})]
\]
and hence
\[
\forall a\in\omega\forall \alpha<a [TI(a,<_{N-1},\Pi_{N}) \to \mbox{{\rm L}}\in M_{2}(M_{2}(\alpha;<_{TW}))]
.\]
\end{theorem}

We see an optimality of this resolving of $\Pi_{N}$-reflecting universes in terms of the lowest recursively Mahlo operation $M_{2}$.

\begin{theorem}\label{th:resolution2}

For each $N\, (2<N<\omega)$ there exist $\Sigma_{1}$ transitive relations $<_{i}\, (2\leq i\leq N-1)$ on 
$\omega$ such that
$<_{N-1}$ is almost wellfounded in {\rm KP}$\ell$, and
{\rm KP}$\Pi_{N}$ is proof-theoretically reducible to
the theory 
\[
\mbox{{\rm KP}}\ell+
\{
\mbox{{\rm L}}\in \bigcap\{M_{2}(M_{2}(\alpha; <_{TW})): 
dom(<_{TW})\ni\alpha<a\}
: a\in\omega\}
\]
for  the restriction $<_{TW}$ of the tower $<_{T}$ of the exponential orderings $<_{E_{i}}$
defined from these to 
the wellfounded parts in the second components hereditarily.
\end{theorem}

Theorem \ref{th:resolution2} is extracted from proof-theoretic analyses of KP$\Pi_{N}$
in \cite{ptpiN} and \cite{WienpiN} .
Let me spend some words on {\it ordinal analyses\/}, an ordinal informative proof-theoretic investigations
in generalities.

\section{Background materials from proof theory}

Let T be a recursive theory containing $\mbox{ACA}_{0}$[the predicative (and hence conservative) extension of the first order arithmetic PA], and $\Pi^{1}_{1}${\it -sound\/}, i.e., 
any T-provable $\Pi^{1}_{1}$-sentence is true in the standard model.

Then its {\it proof-theoretic ordinal\/} $|\mbox{T}|$ is defined to be
the supremum of the order types of the provably recursive well orderings:
\begin{eqnarray*}
&&
|\mbox{T}|  :=   \sup\{\alpha<\omega_{1}^{\footnotesize CK}:  \mbox{T}\vdash WO[<] 
 \,\&\,
 \\
 &&
  \alpha=\mbox{order type } |<| \mbox{ of } < \mbox{ for a recursive ordering } <\}
\end{eqnarray*}
{\bf Remark}.
The ordinal $|\mbox{T}|$ is stable if we consider $\Sigma^{1}_{1}$-orderings and/or 
add true $\Sigma^{1}_{1}$-sentences to $\mbox{T}\supseteq\mbox{ACA}_{0}$,
an anlogue to the C. Spector's boundedness theorem.
For a proof see \cite{attic}.
\\ \smallskip \\
It is seen that $|\mbox{T}|$ is recursive, i.e., $|\mbox{T}|<\omega_{1}^{CK}$, and
easy to cook up a recursive well ordering $<^{T}$ whose order type is equal to $|\mbox{T}|$.

For each $p\in\omega$ let $<_{p}$ denote a recursive well ordering defined as follows:
\begin{enumerate}
\item The case when $p$ is a G\"odel number of a proof in T whose endformula is $WO[\prec]$ for a recursive binary relation $\prec$:
 Then put $<_{p}:=\prec$.
\item Otherwise, let $<_{p}$ denote an empty ordering, i.e., $dom(<_{p})=\emptyset$.
\end{enumerate}
Glue these orderings together to get a recursive ordering $<^{T}$:
\[
\langle n,p\rangle<^{T}\langle m,q\rangle :\Leftrightarrow [p=q \,\&\, n<_{p}m] \lor p<q\]
for a bijective pairing function $\langle n,p\rangle$.

Then $<^{T}$ is a recursive well ordering by the assumptions, and 
\\
$|<^{T}|\leq |T|=\sup\{|<_{p}|:p\in\omega\}\leq|<^{T}|<\omega^{CK}_{1}$ as desired.

Gentzen's celebrated pioneering work yields $|\mbox{ACA}_{0}|=\varepsilon_{0}$.
The first achievement for proof theory of impredicative theory was done by G. Takeuti.
He designed a recursive notation system of ordinals, which describes the proof theoretic ordinal
of, e.g., $\Pi^{1}_{1}$-Comprehension Axiom.
Nowadays Takeuti's proof is understood as for set theories of $\Pi_{2}$-reflecting universes,
i.e., for the Kripke-Platek set theory with the Axiom of Infinity, $\mbox{KP}\omega$.

Ordinal analyses for stronger theories are now obtained.
Let $\langle O(\mbox{T}), <_{T}\rangle$ denote a notation system of proof-theoretic ordinal of 
$\mbox{T}=\mbox{ACA}_{0}$, $\mbox{KP}\omega$, $\mbox{KPM}$, $\mbox{KP}\Pi_{N}$, etc.

Ordinal analyses of theories T show not only the fact $|O(\mbox{{\rm T}})| =|\mbox{{\rm T}}|$ but also more, i.e., some conservative extension results.

\begin{theorem}\label{th:conservation}
Let {\rm EA} denote the elementary recursive arithmetic,
a fragment $I\Delta_{0}+\forall x\exists y(2^{x}=y)$ of {\rm PA}.
\begin{enumerate}
\item\label{th:conservation1}
If $\prec$ is an irreflexive, transitive and provably well founded relation in {\rm T}(not necessarily a total ordering), then there exists an ordinal term $\alpha\in O(\mbox{{\rm T}})$ and an elementary recursive function $f$ so that 
 {\rm EA+}$\forall n,m,k[n\not\prec n \,\&\, (n\prec m\prec k \to n\prec k)]$ proves that 
\[
 \forall n,k[(n\prec k  \to  f(n)<_{T}f(k) )\, \& \, f(n)<_{T}\alpha]
\]

\item\label{th:conservation2}
Over {\rm EA}, $WO[<_{T}]$ is equivalent to the uniform reflection principle 
\\
$\mbox{{\rm RFN}}_{\Pi^{1}_{1}}(\mbox{{\rm T}})$ of {\rm T} for $\Pi^{1}_{1}$-formulas.

\item\label{th:conservation3}
{\rm T} is $\Pi^{1}_{1}$-conservative over the theory 
$\mbox{{\rm ACA}}_{0}\cup\{WO[<_{T}|n] : n\in\omega\}$, which is an extension  of 
$\mbox{{\rm ACA}}_{0}$
by augmenting the wellfoundedness of {\rm each} initial segment $<_{T}|n$ of the ordering $<_{T}$.

\item\label{th:conservation4}
Over {\rm EA},
the 1-consistency $\mbox{{\rm RFN}}_{\Pi^{0}_{2}}(\mbox{{\rm T}})$ of {\rm T} is equivalent to
the fact $\mbox{ERWO}[<_{T}]$ 
that there is no elementary recursive descending chain of ordinals in $O(\mbox{{\rm T}})$.

\item\label{th:conservation5}
{\rm T} is $\Pi^{0}_{2}$-conservative over the theory 
$\mbox{{\rm EA}}\cup\{ERWO[<_{T}|n] : n\in\omega\}$.

Therefore provably recursive functions in {\rm T} are exactly the functions defined by
ordinal recursions along  initial segments $<_{T}|n\, (n\in\omega)$.

\item\label{th:conservation6}
Over {\rm EA},
 finitely iterated consistency statements $\mbox{{\rm CON}}^{(n)}\mbox{{\rm (T)}}$ of {\rm T}
 \[\mbox{{\rm CON}}^{(0)}\mbox{{\rm (T)}} :\Leftrightarrow \forall x(0=0); \: 
\mbox{{\rm CON}}^{(n+1)}\mbox{{\rm (T)}} :\Leftrightarrow \mbox{{\rm CON(T+CON}}^{(n)}\mbox{{\rm (T))}}
\]
 is equivalent to
the inference rule
\[
\infer{A(\alpha)}{[q(\alpha)<_{T}\alpha \to A(q(\alpha))]\to A(\alpha)}
\]
where $\alpha$ denotes a variable ranging over $O(\mbox{{\rm T}})$, and $A$ [$q$] is an elementary recursive relation [function], resp.

\end{enumerate}
\end{theorem}

For a proof of Theorem \ref{th:conservation}.\ref{th:conservation1},
see \cite{attic}.
Theorem \ref{th:conservation}.\ref{th:conservation6} is seen from 
Theorem \ref{th:conservation}.\ref{th:conservation4}
through an Herbrand analysis and a result due to W. Tait\cite{tait}.

The rest of Theorem \ref{th:conservation} is seen from Lemma \ref{lem:conservation} below,
cf. \cite{sod}, \cite{ptrfl}, \cite{ptst}, \cite{ptcoll}, \cite{Leeds}, \cite{odMahlo}, \cite{odpi3}, \cite{ptMahlo}, \cite{ptpi3},
\cite{Wienpi3d} , \cite{ptpiN} and \cite{WienpiN}. 
Also cf. \cite{esubjh}, \cite{esubid}, \cite{esubidea}, \cite{esubMahlo} and \cite{esubpi2} 
for proof theory based on epsilon substitution method.

\begin{lemma}\label{lem:conservation}
\begin{enumerate}
\item\label{lem:conservation1}
{\rm T} proves that each initial segment $<_{T}|n$ is wellfounded.
The proof is uniform in the sense that
\[
\mbox{{\rm EA}}\vdash \mbox{{\rm Proof}}_{T}(p(x),WO[<_{T}|x])
\]
for an elementary recursive function $p(x)$ and a canonical proof predicate $\mbox{{\rm Proof}}_{T}(x,y)$
(read: $x$ is a (code of a) {\rm T}-proof of a (code of a) formula $y$).

\item\label{lem:conservation2}
We can define a rewrite rule(cut-elimination step) $r(p,n)$ on (finite) {\rm T}-proofs $p$ of $\Pi^{1}_{1}$-formulas, 
and an ordinal assignment 
$o: p\mapsto o(p)\in O(\mbox{{\rm T}})$ so that {\rm EA} proves
\[
\forall n[o(r(p,n))<_{T}o(p) \to \mbox{{\rm Tr}}_{\Pi^{1}_{1}}(end(r(p,n)))] \to \mbox{{\rm Tr}}_{\Pi^{1}_{1}}(end(p))
\]
where
$\mbox{{\rm Tr}}_{\Pi^{1}_{1}}$ denotes a partial truth definition for $\Pi^{1}_{1}$-sentences,
and $end(p)$ the end-formula of a proof $p$.

For proofs $p$ of $\Sigma^{0}_{1}$-sentences, the rewrite rule degenerates to be unary, $r(p,n)=r(p,m)$.
\end{enumerate}
\end{lemma}
{\bf NB}.

The size of proof-theoretic ordinals is by no means related to consistency strengths of theories.
Only when we restrict to initial segments of notation systems $O(\mbox{T})$,
the sizes are relevant.
Cf. \cite{lev} and \cite{arnold} for some pathological examples on provably well orderings.

Let 
 $\mbox{CON(T},n\mbox{)}:\Leftrightarrow \forall x\leq n\lnot \mbox{Proof}_{T}(x,\lceil 0=1\rceil)$ denote a partial consistency of T up to $n$.
\begin{enumerate}
\item (\cite{Kreisel})

Let $n\prec m$ denote a recursive relation defined as follows:
\[
n\prec m :\Leftrightarrow
[
\mbox{CON(T},\min\{n,m\}\mbox{)} \,\&\, n<m] \lor 
[
\lnot\mbox{CON(T},\min\{n,m\}\mbox{)} \,\&\, n>m
]
.\]

Even though $|\prec|=\omega$ since T is assumed to be consistent,
$WO[\prec]$ implies CON(T) finitistically.

\item
Modifying the above Kreisel's pathological example, one sees that for any recursive and 
$Bool(\Pi^{1}_{1})$-sound theory T ($Bool(\Pi^{1}_{1})$ denotes the Boolean combinations of $\Pi^{1}_{1}$-sentences), there exists a recursive and 
$Bool(\Pi^{1}_{1})$-sound theory $\mbox{T}^{\prime}$ such that
$|\mbox{T}|<|\mbox{T}^{\prime}|$ but $\mbox{T}^{\prime}\not\vdash \mbox{CON(T)}$:
let $<_{T}$ be any recursive well ordering of type $|\mbox{T}|$, and let
\[
n\prec^{\prime}m :\Leftrightarrow \mbox{CON(T},\max\{n,m\}\mbox{)} \,\&\, n<_{T}m
.\]
Although $|\prec^{\prime}|=|<_{T}|$, $\prec^{\prime}$ is a finite ordering if T is inconsistent.
A fortiori $\mbox{EA}\vdash \lnot\mbox{CON(T)} \to WO[\prec^{\prime}]$.
Hence $\mbox{T}\not\vdash WO[\prec^{\prime}] \to \mbox{CON(T)}$ by
the second incompleteness theorem.
Therefore $\mbox{T}^{\prime}:=\mbox{T}\cup\{WO[\prec^{\prime}]\}$ is a desired one.

Note that if each initial segment of $<_{T}$ is provably wellfounded in T,
then so is for $\prec^{\prime}$.
\end{enumerate}

\section{Collapsing functions iterated}\label{sec:col}

The essential step in cut-elimination for a set theory T is to analyse the axiom expressing
an ordinal $\sigma$ reflects any $\Pi_{2}$-formula $\varphi$:
\[
\varphi^{\mbox{L}_{\sigma}}(a) \land a\in \mbox{L}_{\sigma} \to \exists\beta<\sigma[\varphi^{\mbox{L}_{\beta}} \land a\in \mbox{L}_{\beta}]
.\]
This means that given a proof figure $P$ of the premise, we have to find an ordinal term $\beta<\sigma$:
\[
\begin{array}{ccc}
\infer*[P]{\varphi^{\mbox{L}_{\sigma}}(a) \land a\in \mbox{L}_{\sigma}}{}
&
\Longrightarrow 
&
\infer*{\varphi^{\mbox{L}_{\beta}}(a) \land a\in \mbox{L}_{\beta}}{}
\end{array}
\]

This is done by putting $\beta=d_{\sigma}\alpha<\sigma\, (o(P)=\alpha\in Od(\mbox{T}))$ for a (Mostowski) collapsing function $d$.

Let $C(\alpha)\, (\alpha=o(P))$ denote the set of ordinals which may occur in the reducts of $P$.
Ordinals in $C(\alpha)$ are on the solid lines with gaps here and there
in the following figure:

\vspace{-3cm}
\begin{figure}[h]

\unitlength 3mm
\begin{picture}(12,12)
\put(0,2){\makebox(0,0){$0$}}
\put(0,0){\makebox(0,0){$[$}}
\put(0,0){\thicklines\line(1,0){10}}
\put(10,2){\makebox(0,0){$d_{\sigma}\alpha$}}
\put(10,0){\makebox(0,0){$)$}}
\put(20,2){\makebox(0,0){$\sigma$}}
\put(20,0){\makebox(0,0){$[$}}
\put(20,0){\thicklines\line(1,0){10}}
\put(30,2){\makebox(0,0){$\sigma+d_{\sigma}\alpha$}}
\put(30,0){\makebox(0,0){$)$}}
\put(40,0){\makebox(0,0){$\ldots\ldots\ldots$}}
\end{picture}

\end{figure}

By stuffing the gap below $\sigma$ in the set $C(\alpha)$ up,
$\sigma$ is collapsed down to the least indescribable ordinal $d_{\sigma}\alpha$.
Then ordinals in $C(\alpha)$ cannot discriminate between $\sigma$ and $d_{\sigma}\alpha$
\[\gamma<\sigma\Leftrightarrow \gamma<d_{\sigma}\alpha\, (\gamma\in C(\alpha)),\]
Thus
the ordinal $\beta=d_{\sigma}\alpha$ can be a substitute for  $\sigma$.

To analyse larger ordinals, e.g., $\Pi_{3}$-reflecting ordinals, the collapsing process has to be iterated.

A $\Pi_{3}$-reflecting ordinal $K$ is understood to be $<\varepsilon_{K+1}$-recursively Mahlo, 
$\mbox{L}_{K}\in \bigcap_{\mu<\varepsilon_{K+1}}M^{\mu}_{2}$.
First $K$ is collapsed to a $\mu_{0}$-recursively Mahlo ordinal for a $\mu_{0}<\varepsilon_{K+1}$:
$\kappa_{1}=d^{\mu_{0}}_{K}\alpha_{0}<K$. 
Then $\mbox{L}_{\kappa_{1}}\in M^{\mu_{0}}_{2}$ is collapsed to a $\mu_{1}$-recursively Mahlo ordinal:
$\kappa_{2}=d^{\mu_{1}}_{\kappa_{1}}\alpha_{1}<\kappa_{1}\, (\mu_{1}<\mu_{0})$, etc.
In this way a possibly infinite collapsing process is generated:
$K=\kappa_{0}>d^{\mu_{0}}_{K}\alpha_{0}=\kappa_{1}>d^{\mu_{1}}_{\kappa_{1}}\alpha_{1}=\kappa_{2}>\cdots\,
 (\varepsilon_{K+1}>\mu_{0}>\mu_{1}>\cdots)$.

We have designed a recursive notation system $\langle Od(\Pi_{N}),<\rangle$ of ordinals 
for proof theoretical analysis of KP$\Pi_{N}$, and showed in \cite{ptpiN} that
KP$\Pi_{N}$ is proof-theoretically reducible to the theory
$\mbox{ACA}_{0}+\{WO[<|\alpha]: \Omega>\alpha\in Od(\Pi_{N})\}$,
where $\Omega\in Od(\Pi_{N})$ denotes the least $\Pi_{2}$-reflecting ordinal $\omega_{1}^{\footnotesize CK}$
and $<|\alpha$ the restriction of the ordering $<$ in $Od(\Pi_{N})$ to $\alpha$.
Thus $O(\mbox{KP}\Pi_{N})=Od(\Pi_{N})|\Omega$.

On the other side in \cite{WienpiN} we have shown that
KP$\Pi_{N}$ proves $WO[<|\alpha]$ for {\it each\/} $\alpha<\Omega$.
Indeed, this wellfoundedness proof is essentially formalizable in a theory
$\mbox{{\rm KP}}\ell+
\{
\mbox{{\rm L}}\in \bigcap\{M_{2}(M_{2}(\alpha; <_{TW})): 
dom(<_{TW})\ni\alpha<a\}
: a\in\omega\}$ for some $\Sigma_{1}$ relations $<_{i}\, (2\leq i\leq N-1)$ on $\omega$ such that
$<_{N-1}$ is almost wellfounded in {\rm KP}$\ell$.
This shows Theorem \ref{th:resolution2}.

In the next section we give a sketch of the wellfoundedness proof.

\section{Wellfoundedness proof}

Our wellfoundedness proof of $Od(\Pi_{N})$ is based on the {\it maximal distinguished class\/} 
$\mathcal W$ \cite{Buchholz75}, a $\Sigma_{1}$-definable set of integers, 
and a proper {\it class\/} in KP$\Pi_{N}$. 

To formalize the proof {\it in\/} KP$\Pi_{N}$,
we have to show for each $\eta\in Od(\Pi_{N})$ there exists an $\eta${\it -Mahlo set\/}
on which
the maximal distinguished class enjoys the same closure properties
 as $\mathcal W$ up to the given $\eta$. 
The $\eta$-Mahlo sets are defined through a ramification process 
to resolve the reflecting universes in terms of iterations of lower Mahlo operations\cite{WienpiN}.

\subsection{The notation system $Od(\Pi_{N})$}

The notation system $Od(\Pi_{N})$ (an element of $Od(\Pi_{N})$ is called an {\it ordinal diagram\/},
which is abbreviated o.d.)
contains the constants $\Omega$ for $\omega_{1}^{\footnotesize CK}$ and $\pi$ for the least $\Pi_{N}$-reflecting ordinal.

The main constructor is to form an o.d. $d_{\sigma}^{q}\alpha<\sigma$ from a symbol $d$ and o.d.'s $\sigma,q,\alpha$,
where $\sigma$ denotes a recursively regular ordinal and $q$ a finite sequence of o.d.'s. 

$\gamma\prec_{2}\sigma$ denotes the transitive closure of $\{(\beta,\sigma): \exists\alpha, q(\beta=d_{\sigma}^{q}\alpha)\}$. 
The set $\{\tau:\sigma\prec_{2}\tau\}$ is finite and linearly ordered by $\prec_{2}$ for each $\sigma$, 
namely $\{\sigma:\sigma\preceq_{2}\pi\}$ is a tree with its root $\pi$.

In the diagram $d_{\sigma}^{q}\alpha$, $q$ includes some data telling us how the diagram $d_{\sigma}^{q}\alpha$ is constructed from 
$\{\tau:d_{\sigma}^{q}\alpha\prec_{2}\tau\}=\{\tau:\sigma\preceq_{2}\tau\}$.

The main task in wellfoundedness proofs is to show the tree $\{\sigma:\sigma\preceq_{2}\pi\}$ to be wellfounded. 

Specifically $q$ in $\eta=d_{\sigma}^{q}\alpha$ includes some data $st_{i}(\eta), pd_{i}(\eta), rg_{i}(\eta)$ for $2\leq i<N$. 
$st_{N-1}(\eta)$ is an o.d. less than $\varepsilon_{\pi+1}$, and $pd_{2}(\eta)=\sigma$.

A relation $\prec_{i}$ is defined from $pd_{i}(\eta)$ as the transitive closure of 
$\{(\eta,\kappa):\kappa=pd_{i}(\eta)\}$.
This enjoys 
$\prec_{i+1}\subseteq\prec_{i}$.
Therefore the diagram $pd_{i}(\eta)$ is a proper subdiagram of $\eta$. 
 $st_{i}(\eta)$ is an o.d. less than the next admissible $\kappa^{+}$ to a $\kappa=rg_{i}(\eta)\leq pd_{i+1}(\eta)$.
$rg_{N-1}(\eta)=\pi$ for any such $\eta=d_{\sigma}^{q}\alpha$.

$q$ determines a sequence $\{\eta_{i}^{m}: m<lh_{i}(\eta)\}$ of subdiagrams of $\eta$
 with its length $lh_{i}(\eta)=n+1>0$. The sequence enjoys the following property:
\[
\eta\preceq_{i+1}\eta_{i}^{0}\prec_{i+1}\eta_{i}^{1}\prec_{i+1}\cdots\prec_{i+1}\eta_{i}^{n}<\pi
\]
with $st_{i}(\eta_{i}^{m})<(rg_{i}(\eta_{i}^{m}))^{+}$.

\subsection{Towers derived from ordinal diagrams}

Define relations $\ll_{i}$ for $2\leq i\leq N-1$ by
\[
\eta\ll_{i}\rho:\Leftrightarrow \eta\prec_{i}\rho \,\&\, rg_{i}(\eta)=rg_{i}(\rho) \,\&\,
st_{i}(\eta)<st_{i}(\rho)
.\]
Extend $\ll_{i}$ by augmenting the least element $1$:
\[
1\ll_{i}\eta
.\]
$\pi^{\alpha}$ denotes $\pi^{\alpha}\cdot 1$.

Let $\lhd_{i}:\equiv <_{E_{i}}$ be exponential ordering defined from $\ll_{i}\, (2\leq i\leq N-1)$.
Namely $\lhd_{N-1}:\equiv\ll_{N-1}$ and $\lhd_{i}:\equiv E(\lhd_{i+1},\ll_{i})$, cf. Definition \ref{df:LEe}.

Extend $\lhd_{i}$ to $\lhd_{i}^{+}$ by adding the successor function $+1$.
Namely the domain is expanded to $dom(\lhd_{i}^{+}):=dom(\lhd_{i})\cup\{a+1: a\in dom(\lhd_{i})\}$,
and define for $a,b\in dom(\lhd_{i})$
\begin{eqnarray*}
a+1\lhd_{i}^{+}b+1 & :\Leftrightarrow & a\lhd_{i}b \\
a+1\lhd_{i}^{+}b & :\Leftrightarrow & a\lhd_{i}b \\
a\lhd_{i}^{+}b+1 & :\Leftrightarrow & a\lhd_{i}b \mbox{ or } a=b
\end{eqnarray*}

From the sequence $\{\eta_{i}^{m}: 2\leq i<N-1, m<lh_{i}(\eta)\}$ we define a tower $T(\eta)=E_{2}(\eta)$.
The elements of the form $E_{i}(\eta)(+1)$ are understood to be ordered by $\lhd^{+}_{i}$.
Let $\lhd_{T}:\equiv\lhd^{+}_{2}$.

\begin{eqnarray*}
E_{N-1}(\eta) & := & \eta \\
E_{i}(\eta) & := & \sum_{1\leq m< lh_{i}(\eta)}\pi^{E_{i+1}(\eta_{i}^{m})} \eta_{i}^{m-1}
+\pi^{E_{i+1}(\eta_{i}^{0})+1}+\pi^{E_{i+1}(\eta)}
\end{eqnarray*}

The sequence $\{\eta_{i}^{m}: m<lh_{i}(\eta)\}$ is defined so that, cf. \cite{WienpiN} for a proof,
\[
\gamma\prec_{i}\eta \Rightarrow E_{i}(\gamma)\lhd^{+}_{i}E_{i}(\eta)
.\]
In particular
\begin{equation}\label{eq:tower2}
\gamma\prec_{2}\eta \Rightarrow T(\gamma)\lhd_{T} T(\eta)
\end{equation}

\subsection{Distinguished classes}

An elementary fact on the maximal distinguished class $\mathcal W$ says that 
$\mathcal W$ is well ordered by $<$ on 
$Od(\Pi_{N})$, and
$\mathcal W|\Omega$ is included in the wellfounded part of $Od(\Pi_{N})$.
Therefore it suffices to  show  $\eta\in\mathcal W$ for {\it each\/} $\eta\in Od(\Pi_{N})$.

$\mathcal W$ is defined to be the union of the distinguished sets, 
\[
\mathcal W=\bigcup\{X\subseteq Od(\mbox{T}): D[X]\}
\]
where $D[X]$(read:$X$ is a distinguished set) is a $\Delta_{1}$-formula on limits of admissible sets.
Hence $\mathcal W$ is a $\Sigma_{1}$-definable set of integers, and a proper {\it class\/} in 
KP$\Pi_{N}$.

Since $D[X]$ is $\Delta_{1}$ on limits of admissibles, it is absolute: 
$D[X] \Leftrightarrow P\models D[X]$ for any $X\in P\cap\mathcal{P}(\omega)$.
Let 
$\mathcal W^{P}=\bigcup\{X\in P: P\models D[X]\}$
 denote the maximal distinguished class on $P$.

The following is a key on distinguished sets.

\begin{lemma}\label{lem:3wf6}

There exists a $\Pi_{2}$-formula $g(\eta)\, (\eta\in Od(\Pi_{N}))$ for which the following holds
for any limits $Q$ of admissibles:
Assume $g(\eta)^{Q}$ and
\begin{equation}\label{eq:etamahlo}
\forall\gamma\prec_{2}\eta\{g(\gamma)^{Q}  \Rightarrow   \gamma\in \mathcal W^{Q}\}
\end{equation}
Then there exists a distinguished class $X$ such that $\eta\in X$ and $X$ is definable in $Q$.
\end{lemma}

For some $\Sigma_{1}$ classes $U_{i}$ on $\omega$, 
the $\Sigma_{1}$ transitive relations on $\omega$, $<_{i}$ mentioned in Theorem \ref{th:resolution2}
are now defined to be 
\[
\eta<_{i}\rho:\Leftrightarrow \eta\ll_{i}\rho \,\&\, \eta,\rho\in U_{i}
.\]
By definition $1\in U_{i}$ for any $i$.
$<_{N-1}$ is seen to be almost wellfounded in KP$\ell$.
 
Let $<_{TW}$ denote  the restriction of the tower $<_{T}$ of the exponential orderings $<_{E_{i}}$
defined from these $\Sigma_{1}$ relations $<_{i}\, (2\leq i\leq N-1)$ to 
the wellfounded parts in the second components hereditarily.

In other words, 
\[
T(\eta)<_{T}T(\rho) \Leftrightarrow T(\eta)\lhd_{T}T(\rho) \,\&\, 
\forall i[{\cal K}_{i}(\eta)\cup{\cal K}_{i}(\rho)\subseteq U_{i}]
\]
and
\[
T(\eta)<_{TW}T(\rho) \Leftrightarrow T(\eta)<_{T}T(\rho) \,\&\,
 \forall i>0[{\cal K}_{i}(\eta)\cup{\cal K}_{i}(\rho)\subseteq W(<_{i})]
 \]
 where
 \begin{enumerate}
 \item
 ${\cal K}_{2}(\eta):=\{\eta_{2}^{m}: m<lh_{2}(\eta)\}$.
 \item
 For $2<i<N-1$,
 ${\cal K}_{i}(\eta):=\{\rho_{i}^{m}: m<lh_{i}(\rho), \rho\in {\cal K}_{i-1}(\eta)\}$.
\end{enumerate}

\begin{lemma}\label{th:5awf16.etamahlo}
If $P\in M_{2}(M_{2}(T(\eta);<_{TW}))$, then $g(\eta)^{P} \to \eta\in\mathcal W^{P}$.
\end{lemma}
{\bf Proof} by induction on $\in$.
Suppose $P\in M_{2}(M_{2}(T(\eta);<_{TW}))$ and $g(\eta)^{P}$.
Pick a $Q\in P\cap M_{2}(T(\eta);<_{TW})$ so that $g(\eta)^{Q}$.

We show (\ref{eq:etamahlo}).
Assume $\gamma\prec_{2}\eta$ and $g(\gamma)^{Q}$.
(\ref{eq:tower2}) yields $T(\gamma)\lhd_{T}T(\eta)$.
On the other side the $\Pi_{2}$ formula $g(\gamma)$ is defined so that
\[
g(\gamma)^{Q} \to \forall i[{\cal K}_{i}(\gamma)\subseteq U^{Q}_{i}] \,\&\, \forall i>0[{\cal K}_{i}(\gamma)\subseteq W^{Q}(<^{Q}_{i})]
.\]
Since $\bigcup_{i}{\cal K}_{i}(\eta)$ is finite, we can assume $\forall i[{\cal K}_{i}(\eta)\subseteq U^{Q}_{i}]$,
and hence
$T(\gamma)<_{TW}^{Q}T(\eta)$.
Therefore $Q\in M_{2}(M_{2}(T(\gamma);<_{TW}))$.
IH yields $\gamma\in \mathcal W^{Q}$. 
This shows (\ref{eq:etamahlo}).

By Lemma \ref{lem:3wf6}, let $X$ be a distinguished class
 definable over $Q$ such that $\eta\in X$.
Thus $X\in P \,\&\, D[X]$, and $\eta\in \mathcal W^{P}$.
\hspace*{\fill} $\Box$

Assuming $\mbox{L}\in M_{2}(M_{2}(T(\eta);<_{TW}))$ for each $\eta$,
we have $g(\eta)^{\mbox{\footnotesize L}} \to \eta\in\mathcal W^{\mbox{\footnotesize L}}=\mathcal W$
 by Lemma \ref{th:5awf16.etamahlo}.
On the other side, it is not hard to show $g(\eta)^{\mbox{\footnotesize L}}$ for  each $\eta$
 in KP$\ell$.
 
Therefore the wellfoundedness of $Od(\Pi_{N})$ up to each $\eta<\Omega$
follows from 
$\{
\mbox{{\rm L}}\in M_{2}(M_{2}(T(\eta); <_{TW}))
: \eta\in Od(\Pi_{N})\}$ over KP$\ell$.


\begin{thebibliography}{99}

\bibitem[A96a]{sod}T. Arai, Systems of ordinal diagrams, draft, 1996.
\bibitem[A96b]{ptrfl} T. Arai, Proof theory for theories of ordinals I: Reflecting ordinals, draft, 1996.
\bibitem[A97a]{ptst} T. Arai, Proof theory for theories of ordinals II: $\Sigma_{1}$-stability, draft, 1997.
\bibitem[A97b]{ptcoll} T. Arai, Proof theory for theories of ordinals III: $\Pi_{1}$-collection, draft, 1997.

\bibitem[A98]{attic} T. Arai,  Some results on cut-elimination, provable well-orderings, induction and reflection, 
Ann. Pure Appl. Logic 95 (1998) 93-184.

\bibitem[A99]{Leeds} T. Arai,  Introduction to finitary analyses of proof figures, 
In: Sets and Proofs. Invited papers from Logic Colloquium '97-European Meeting of the Association for Symbolic Logic, Leeds, July 1997. Ed. by S. B. Cooper and J. K. Truss,
London Mathematical Society Lecture Notes, vol. 258, Cambridge University Press (1999), pp.1-25.

\bibitem[A00a]{odMahlo} T. Arai, Ordinal diagrams for recursively Mahlo universes, 
Arch. Math. Logic 39 (2000) 353-391.

\bibitem[A00b]{odpi3} T. Arai, Ordinal diagrams for $\Pi_{3}$-reflection, 
Jour. Symb. Logic 65 (2000) 1375-1394.


\bibitem[A02]{esubjh} T. Arai, Epsilon substitution method for theories of jump hierarchies, 
 Arch. Math. Logic 41 (2002) 123-153.
 
\bibitem[A03a]{esubid} T. Arai, Epsilon substitution method for $ID_{1}(\Pi^{0}_{1}\lor\Sigma_{1}^{0})$, 
Ann. Pure Appl. Logic 121 (2003) 163-208.

\bibitem[A03b]{ptMahlo} T. Arai, Proof theory for theories of ordinals I:recursively Mahlo ordinals, 
Ann. Pure Appl. Logic 122 (2003) 1-85.

\bibitem[A04a]{ptpi3} T. Arai, Proof theory for theories of ordinals II:$\Pi_{3}$-Reflection, 
Ann. Pure Appl. Logic 129 (2004) 39-92.

\bibitem[A04b]{Wienpi3d} T. Arai, Wellfoundedness proofs by means of non-monotonic inductive definitions I: 
$\Pi^{0}_{2}$-operators, 
Jour. Symb. Logic 69 (2004) 830-850.

\bibitem[A05a]{esubidea}T. Arai,  Ideas in the epsilon substitution method for $\Pi^{0}_{1}$-FIX,
Ann. Pure Appl. Logic 136 (2005) 3-21.

\bibitem[A05b]{esubMahlo} T. Arai, Epsilon substitution method for $[\Pi^{0}_{1},\Pi^{0}_{1}]$-FIX, 
Arch. Math. Logic 44 (2005) 1009-1043.

\bibitem[A06]{esubpi2} T. Arai, Epsilon substitution method for $\Pi^{0}_{2}$-FIX, 
Jour. Symb. Logic 71 (2006) 1155-1188.


\bibitem[A$\infty$a]{ptpiN} T. Arai, Proof theory for theories of ordinals III:$\Pi_{N}$-reflection, submitted.

\bibitem[A$\infty$b]{WienpiN} T. Arai, Wellfoundedness proofs by means of non-monotonic inductive definitions II: 
first order operators, submitted.



\bibitem[Beckmann02]{arnold} A. Beckmann, A non-well-founded primitive recursive tree provably well-founded for co-r.e. sets, Arch. Math. Logic 41(2002) 251-257.

\bibitem[Beklemishev00]{lev} L. Beklemishev, Another pathological well-ordering, in Logic Colloquium 98(Prague), 105-108, Lect. Notes Logic 13, Assoc. Symb. Logic, 2000.

\bibitem[Buchholz75]{Buchholz75} W. Buchholz, Normalfunktionen und konstruktive Systeme von Ordinalzahlen. In: Diller, J., M\"uller, G.H.(eds.) Proof Theory Symposion, Kiel 1974 (Lecture Notes in Mathematics, vol.500, pp.4-25). Berlin: Springer 1975

\bibitem[Kreisel77]{Kreisel} G. Kreisel, Wie die Beweistheoire zu ihren Ordinalzahlen kam und kommt, 
Jber. Deutsch. Math.-Verein 78(1977), 177-223.

\bibitem[Richter-Aczel74]{Richter-Aczel74} W.H. Richter and  P. Aczel, Inductive definitions and reflecting properties of admissible ordinals, Generalized Recursion Theory, Studies in Logic, vol.79, North-Holland, 1974, pp.301-381.

\bibitem[Tait65]{tait} W.  W. Tait, Functionals defined by transfinite recursion, 
Jour. Symb. Logic  30 (1965) 155-174.
\end{thebibliography}
\end{document}